\documentclass[12pt]{article}
\usepackage{mathrsfs}
\usepackage{amssymb}
\usepackage{amsmath}
\usepackage{amsbsy}
\usepackage{epsfig}
\usepackage{enumerate}
\usepackage{bm}
\newtheorem{lemma}{Lemma}[section]
\newtheorem{theorem}{Theorem}[section]
\newtheorem{proposition}{Proposition}[section]

\newtheorem{definition}{Deriniton}[section]
\newtheorem{remark}{Remark}[section]

\newbox\TempBox \newbox\TempBoxA

\def\ep{\textsf{E}} 
\def\Sbep{\widehat{\mathbb E}} 
\def\cSbep{\widehat{\mathcal E}} 
\def\Capc{\mathbb V} 
\def\cCapc{\mathcal V} 

\def\underwiggle 1{
\ifmmode\setbox\TempBox=\hbox{$ 1$}\else\setbox\TempBox=\hbox{
1}\fi \setbox\TempBoxA=\hbox to \wd\TempBox{\hss\char'176\hss}
\rlap{\copy\TempBox}\smash{\lower9pt\hbox{\copy\TempBoxA}} }
\renewcommand{\baselinestretch}{1.5}

\begin{document}

\thispagestyle{empty}

\begin{center}
 { \LARGE\bf Convergence to a self-normalized G-Brownian motion$^{\ast}$}
\end{center}

\begin{center} {\sc
Li-Xin Zhang\footnote{Research supported by Grants from the National Natural Science Foundation of China (No. 11225104)  and the Fundamental Research Funds for the Central Universities.
}
}\\
{\sl \small Department of Mathematics, Zhejiang University, Hangzhou 310027} \\
(Email:stazlx@zju.edu.cn)\\
\today
\end{center}

\renewcommand{\abstractname}{~}
\begin{abstract}
\centerline{\bf Abstract}
G-Brownian motion has a very rich and interesting new structure which nontrivially
generalizes the classical one. Its
quadratic variation process is also a continuous process with independent
and stationary increments. We prove a self-normalized functional central limit theorem  for independent and identically distributed random variables under the sub-linear expectation with the  limit process being a  G-Browian motion self-normalized by its quadratic variation.
To prove the self-normalized central limit theorem, we also establish  a new  Donsker's invariance principle with the  limit process being a  generalized G-Browian motion.


\end{abstract}

\baselineskip 22pt

\renewcommand{\baselinestretch}{1.5}



\section{ Introduction}\label{sectIntro}
\setcounter{equation}{0}

Let $\{X_n; n\ge 1\}$ be a sequence of independent and identically distributed random variables on a probability space $(\Omega, \mathscr{F}, P)$.  Set $S_n=\sum_{j=1}^n X_j$. Suppose   $EX_1=0$ and $EX_1^2=\sigma^2>0$. The well-known central limit theorem says that
\begin{equation}\label{classicalCLT1} \frac{S_n}{\sqrt{n}}\overset{d}\to N(0,\sigma^2),
\end{equation}
or equivalently, for any bounded continuous function $\psi(x)$,
\begin{equation}\label{classicalCLT2} E\left[\psi\Big( \frac{S_n}{\sqrt{n}}\Big)\right]\to E\left[\psi(\xi)\right],
\end{equation}
where $ \xi\sim N(0,\sigma^2)$ is a normal random variable. If the normalization factor $\sqrt{n}$ is replaced by $\sqrt{V_n}$ where $V_n=\sum_{j=1}^n X_j^2$,  then
\begin{equation}\label{classicalCLT3} \frac{S_n}{\sqrt{V_n}}\overset{d}\to N(0,1).
\end{equation}
Gin\'e, G\"otze  and Mason (1997)   proved that (\ref{classicalCLT3}) holds if and only if $EX_1=0$ and
\begin{equation}\label{classicalCLT4} \lim_{x\to \infty} \frac{x^2P\big(|X_1|\ge x\big)}{EX_1^2I\{|X_1|\le x\}}=0.
 \end{equation}
 The result (\ref{classicalCLT3}) is  refereed to as the self-normalized central limit theorem. The purpose of this paper is to establish the self-normalized central limit theorem under the sub-linear expectation.

The sub-linear expectation  or called  and G-expectation is a nonlinear expectation advancing
the notions of  backward stochastic differential equations,  g-expectations, and provides  a   flexible framework to model non-additive probability  problems  and the volatility uncertainty in finance.
Peng (2006, 2008a, 2008b) introduced a general framework of the sub-linear expectation  of random variables    and introduced the notions of  G-normal random variable, G-Brownian motion,  independent and identically distributed random variables etc under  the sub-linear expectations. The construction of sub-linear expectations on   the space of continuous
paths and discrete time paths can also be founded in   Yan et al (2012) and   Nutz and  Handel (2013).  For basic  properties of the sub-linear expectations, one can refer to Peng (2008b, 2009, 2010a, etc). For  stochastic
calculus and stochastic differential equations with respect to a G-Brownian motion, one can refer  Li and Peng (2011), Hu, et al (2014a,b) etc and a book of Peng (2010a).

  The central limit theorem under the sub-linear expectation was first established by Peng (2008b). It says that (\ref{classicalCLT2}) remains true when the expectation $E$ is replaced by a sub-linear expectation $\Sbep$ if  $\{X_n;n\ge 1\}$ are   independent and identically distributed   under $\Sbep$, i.e.,
\begin{equation}\label{CLT1} \frac{S_n}{\sqrt{n}}\overset{d}\to \xi \text{ under } \Sbep,
\end{equation}
  where  $\xi$ is    a G-normal random variable.

  In the classical   case, when $\ep[X_1^2]$ is finite, (\ref{classicalCLT3}) follows from the cental limit theorem (\ref{classicalCLT1}) immediately by Slutsky's lemma and the fact that
  $$ \frac{V_n}{n}\overset{P}\to \sigma^2. $$
  The later is due to the law of large numbers.  In the framework of the sub-linear expectation,  $\frac{V_n}{n}$ no longer converges to a constant.
  The self-normalized central limit theorem can not follow from the central limit theorem (\ref{CLT1}) directly.
In this paper,  we will prove that
\begin{equation}\label{CLT2} \frac{S_n}{\sqrt{V_n}}\overset{d}\to \frac{W_1}{\sqrt{\langle W\rangle_1}} \text{ under } \Sbep,
\end{equation}
where $W_t$ is a G-Brownian motion and $\langle W \rangle_t$ is its quadratic variation process. A very interesting  phenomenon of G-Brownian motion is that its
quadratic variation process is also a continuous process with independent
and stationary increments, and thus can be still regarded as a Brownian motion. When the sub-linear expectation $\Sbep$ reduces to a linear one,   $W_t$ is the classical Borwnian motion with $W_1\sim N(0,\sigma^2)$ and $\langle W \rangle_t=t\sigma^2$, and then  (\ref{CLT2}) is just (\ref{classicalCLT3}).  Our  main results on the self-normalized central limit theorem   will be given in Section \ref{sectMain}  where the process of the self-normalized partial sums ${S_{[nt]}}/{\sqrt{V_n}}$ is proved to converge to a self-normalized G-Browian motion ${W_t}/{\sqrt{\langle W\rangle_1}}$. We also consider the case that the second moments of $X_i$s are infinite and obtain  the  self-normalized central limit theorem under a condition similar to (\ref{classicalCLT4}).
In the next section, we state basic settings in a sub-linear expectation space including, capacity, independence, identical distribution, G-Brownian motion etc. One can skip this section if he/she is familiar with these  concepts.  To prove the self-normalized central limit theorem, we establish  a new  Donsker's invariance principle in Section \ref{sectIP} with the limit process being a generalized G-Browian motion. The proof is given in the last section.

\section{Basic Settings}\label{sectBasic}
\setcounter{equation}{0}

We use the framework and notations of Peng (2008b). Let  $(\Omega,\mathcal F)$
 be a given measurable space  and let $\mathscr{H}$ be a linear space of real functions
defined on $(\Omega,\mathcal F)$ such that if $X_1,\ldots, X_n \in \mathscr{H}$  then $\varphi(X_1,\ldots,X_n)\in \mathscr{H}$ for each
$\varphi\in C_b(\mathbb R^n)\bigcup  C_{l,Lip}(\mathbb R^n)$,  where $C_b(\mathbb R^n)$ denote  the
space  of all bounded  continuous functions and $C_{l,Lip}(\mathbb R^n)$ denotes the linear space of (local Lipschitz)
functions $\varphi$ satisfying
\begin{eqnarray*} & |\varphi(\bm x) - \varphi(\bm y)| \le  C(1 + |\bm x|^m + |\bm y|^m)|\bm x- \bm y|, \;\; \forall \bm x, \bm y \in \mathbb R^n,&\\
& \text {for some }  C > 0, m \in \mathbb  N \text{ depending on } \varphi. &
\end{eqnarray*}
$\mathscr{H}$ is considered as a space of ``random variables''. In this case we denote $X\in \mathscr{H}$.
Further, we let $C_{b,Lip}(\mathbb R^n)$  denote  the
space  of all bounded and Lipschitz functions   on $\mathbb R^n$.

\subsection{Sub-linear expectation and capacity}
\begin{definition}\label{def1.1} A {\bf sub-linear expectation} $\Sbep$ on $\mathscr{H}$  is a functional $\Sbep: \mathscr{H}\to \overline{\mathbb R}$ satisfying the following properties: for all $X, Y \in \mathscr H$, we have
\begin{description}
  \item[\rm (a)] {\bf Monotonicity}: If $X \ge  Y$ then $\Sbep [X]\ge \Sbep [Y]$;
\item[\rm (b)] {\bf Constant preserving}: $\Sbep [c] = c$;
\item[\rm (c)] {\bf Sub-additivity}: $\Sbep[X+Y]\le \Sbep [X] +\Sbep [Y ]$ whenever $\Sbep [X] +\Sbep [Y ]$ is not of the form $+\infty-\infty$ or $-\infty+\infty$;
\item[\rm (d)] {\bf Positive homogeneity}: $\Sbep [\lambda X] = \lambda \Sbep  [X]$, $\lambda\ge 0$.
 \end{description}
 Here $\overline{\mathbb R}=[-\infty, \infty]$. The triple $(\Omega, \mathscr{H}, \Sbep)$ is called a sub-linear expectation space. Give a sub-linear expectation $\Sbep $, let us denote the conjugate expectation $\cSbep$of $\Sbep$ by
$$ \cSbep[X]:=-\Sbep[-X], \;\; \forall X\in \mathscr{H}. $$
\end{definition}

Next, we introduce the capacities corresponding to the sub-linear expectations.
Let $\mathcal G\subset\mathcal F$. A function $V:\mathcal G\to [0,1]$ is called a capacity if
$$ V(\emptyset)=0, \;V(\Omega)=1 \; \text{ and } V(A)\le V(B)\;\; \forall\; A\subset B, \; A,B\in \mathcal G. $$
It is called to be sub-additive if $V(A\bigcup B)\le V(A)+V(B)$ for all $A,B\in \mathcal G$  with $A\bigcup B\in \mathcal G$.

 Let $(\Omega, \mathscr{H}, \Sbep)$ be a sub-linear space,
and  $\cSbep $  be  the conjugate expectation of $\Sbep$. It is natural to define the capacity of a set $A$ to be the sub-linear expectation of
the   indicator function $I_A$ of $A$. However, $I_A$ may be not in $\mathscr{H}$. So, we denote a pair $(\Capc,\cCapc)$ of capacities by
$$ \Capc(A):=\inf\{\Sbep[\xi]: I_A\le \xi, \xi\in\mathscr{H}\}, \;\; \cCapc(A):= 1-\Capc(A^c),\;\; \forall A\in \mathcal F, $$
where $A^c$  is the complement set of $A$.
Then $\Capc$  is sub-additive and
\begin{equation}\label{eq1.3} \begin{matrix}
&\Capc(A)=\Sbep[I_A], \;\; \cCapc(A)= \cSbep[I_A],\;\; \text{ if } I_A\in \mathscr H\\
&\Sbep[f]\le \Capc(A)\le \Sbep[g], \;\;\cSbep[f]\le \cCapc(A) \le \cSbep[g],\;\;
\text{ if } f\le I_A\le g, f,g \in \mathscr{H}.
\end{matrix}
\end{equation}

Further, we define an extension of $\Sbep^{\ast}$ of $\Sbep$ by
$$ \Sbep^{\ast}[X]=\inf\{\Sbep[Y]:X\le Y, \; Y\in \mathscr{H}\}, \;\; \forall X:\Omega\to \mathbb R, $$
where $\inf\emptyset=+\infty$. Then
$$ \begin{matrix}
&\Sbep^{\ast}[X]=\Sbep[X]\; \text{ if } X\in \mathscr H, \;\; \;\; \Capc(A)=\Sbep^{\ast}[I_A], \\
&\Sbep[f]\le \Sbep^{\ast}[X]\le \Sbep[g]\;\;
\text{ if } f\le X\le g, f,g \in \mathscr{H}.
\end{matrix}
$$

\subsection{Independence and distribution}
\begin{definition} {\em(Peng (2006, 2008b))}

\begin{description}
  \item[ \rm (i)] ({\bf Identical distribution}) Let $\bm X_1$ and $\bm X_2$ be two $n$-dimensional random vectors defined
respectively in sub-linear expectation spaces $(\Omega_1, \mathscr{H}_1, \Sbep_1)$
  and $(\Omega_2, \mathscr{H}_2, \Sbep_2)$. They are called identically distributed, denoted by $\bm X_1\overset{d}= \bm X_2$  if
$$ \Sbep_1[\varphi(\bm X_1)]=\Sbep_2[\varphi(\bm X_2)], \;\; \forall \varphi\in C_{l,Lip}(\mathbb R^n), $$
whenever the sub-expectations are finite. A sequence $\{X_n;n\ge 1\}$ of random variables is said to be identically distributed if $X_i\overset{d}= X_1$ for each $i\ge 1$.
\item[\rm (ii)] ({\bf Independence})   In a sub-linear expectation space  $(\Omega, \mathscr{H}, \Sbep)$, a random vector $\bm Y =
(Y_1, \ldots, Y_n)$, $Y_i \in \mathscr{H}$ is said to be independent to another random vector $\bm X =
(X_1, \ldots, X_m)$ , $X_i \in \mathscr{H}$ under $\Sbep$  if for each test function $\varphi\in C_{l,Lip}(\mathbb R^m \times \mathbb R^n)$
we have
$$ \Sbep [\varphi(\bm X, \bm Y )] = \Sbep \big[\Sbep[\varphi(\bm x, \bm Y )]\big|_{\bm x=\bm X}\big],$$
whenever $\overline{\varphi}(\bm x):=\Sbep\left[|\varphi(\bm x, \bm Y )|\right]<\infty$ for all $\bm x$ and
 $\Sbep\left[|\overline{\varphi}(\bm X)|\right]<\infty$.
 \item[\rm (iii)] ({\bf IID random variables}) A sequence of random variables $\{X_n; n\ge 1\}$
 is said to be independent and identically distributed (IID), if
 $X_i\overset{d}=X_1$ and $X_{i+1}$ is independent to $(X_1,\ldots, X_i)$ for each $i\ge 1$.
 \end{description}
\end{definition}

\subsection{G-normal distribution, G-Brownian motion and its quadratic variation}
Let $0<\underline{\sigma}\le \overline{\sigma}<\infty$ and $G(\alpha)=\frac{1}{2}(\overline{\sigma}^2 \alpha^+ - \underline{\sigma}^2 \alpha^-)$. $X$ is call a normal $N\big(0, [\underline{\sigma}^2, \overline{\sigma}^2]\big)$ distributed random variable  (write $X\sim N\big(0, [\underline{\sigma}^2, \overline{\sigma}^2]\big)$) under $\Sbep$, if for any bounded Lipschitz function $\varphi$, the function $u(x,t)=\Sbep\left[\varphi\left(x+\sqrt{t} X\right)\right]$ ($x\in \mathbb R, t\ge 0$) is the unique viscosity solution of  the following heat equation:
      $$ \partial_t u -G\left( \partial_{xx}^2 u\right) =0, \;\; u(0,x)=\varphi(x). $$

 Let $C[0,1]$ be a function space of continuous functions   on $[0,1]$ equipped with the super-norm $\|x\|=\sup\limits_{0\le t\le 1}|x(t)|$ and $C_b\big(C[0,1]\big)$ is the set of bounded continuous  functions $h(x):C[0,1]\to \mathbb R$. The modulus of the continuity of an element   $x\in C[0,1]$  is defined by $$\omega_{\delta}(x)=\sup_{|t-s|<\delta}|x(t)-x(s)|.$$

It is showed that  there is a sub-linear expectation space $\big(\widetilde{\Omega}, \widetilde{\mathscr{H}},\widetilde{\mathbb E}\big)$ with
$\widetilde{\Omega}= C[0,1]$ and $C_b\big(C[0,1]\big)\subset \widetilde{\mathscr{H}}$ such that $(\widetilde{\mathscr{H}}, \widetilde{\mathbb E}[\|\cdot\|])$ is a Banach space, and
the canonical process $W(t)(\omega) = \omega_t  (\omega\in \widetilde{\Omega})$ is a G-Brownian motion with
$W(1)\sim N\big(0, [\underline{\sigma}^2, \overline{\sigma}^2]\big)$ under $\widetilde{\mathbb E}$, i.e.,
  for all $0\le t_1<\ldots<t_n\le 1$, $\varphi\in C_{l,lip}(\mathbb R^n)$,
\begin{equation}\label{eqBM} \widetilde{\mathbb E}\left[\varphi\big(W(t_1),\ldots, W(t_{n-1}), W(t_n)-W(t_{n-1})\big)\right]
  =\widetilde{\mathbb E}\left[\psi\big(W(t_1),\ldots, W(t_{n-1})\big)\right],
  \end{equation}
  where $\psi(x_1,\ldots, x_{n-1})\big)=\widetilde{\mathbb E}\left[\varphi\big(x_1,\ldots, x_{n-1}, \sqrt{t_n-t_{n-1}}W(1)\big)\right]$
  (c.f. Peng (2006, 2008a, 2010), Denis, Hu and Peng (2011)).

  The quadratic variation process of a G-Brownian motion $W$ is defined by
  $$\langle W \rangle_t=\lim_{\|\Pi_t^N\|\to 0}\sum_{j=1}^{N-1}\big(W(t_j^N)-W(t_{j-1}^N)\big)^2=W^2(t)-2\int_0^t W(t) dW(t),$$
  where $\Pi_t^N=\{t_0^N,t_1^N,\ldots, t_N^n\}$ is a partition of $[0,t]$ and $\|\Pi_t^N\|=\max_j|t_j^N-t_{j-1}^N|$, and the limit is taken in $L_2$, i.e.,
  $$ \lim_{\|\Pi_t^N\|\to 0}\widetilde{\mathbb E}\left[\Big(\sum_{j=1}^{N-1}\big(W(t_j^N)-W(t_{j-1}^N)\big)^2-\langle W \rangle_t\Big)^2\right]=0. $$
The quadratic variation process $\langle W \rangle_t$ is also a continuous process with independent
and stationary increments. For the properties and the distribution of the quadratic variation process, one can refer to a book of Peng (2010a).

Denis, Hu and Peng (2011)   showed the following representation of the G-Brownian motion (c.f, Theorem 52).
 \begin{lemma} \label{DenisHuPeng}
 Let $(\Omega, \mathscr{F}, P) $ be a probability
measure space and $\{B(t)\}_{t\ge 0}$  is a $P$-Brownian motion. Then for all bounded continuous function $\varphi: C_b[0,1]\to \mathbb R$,
$$ \widetilde{\mathbb E}\left[\varphi\big(W(\cdot)\big)\right]=\sup_{\theta\in \Theta}\ep_P\left[\varphi\big(W_{\theta}(\cdot)\big)\right],\;\;
W_{\theta}(t) = \int_0^t\theta(s) dB(s), $$
where
\begin{eqnarray*}
&\Theta=\left\{ \theta:\theta(t) \text{ is } \mathscr{F}_t\text{-adapted process such that }  \underline{\sigma}\le \theta(t)\le \overline{\sigma}\right\},&\\
& \mathscr{F}_t=\sigma\{B(s):0\le s\le t\}\vee \mathscr{N}, \;\; \mathscr{N} \text{ is the collection of } P\text{-null subsets}. &
\end{eqnarray*}
\end{lemma}

 In the sequel of this paper, the sequences $\{X_n;n\ge 1\}$, $\{Y_n;n\ge 1\}$ etc of the random variables are considered in $(\Omega, \mathscr{H}, \Sbep)$. Without specification, we suppose that $\{X_n; n\ge 1\}$ is   a sequence of independent and identically distributed random variables in $(\Omega, \mathscr{H}, \Sbep)$  with  $\Sbep[X_1]=\cSbep[X_1]=0$, $\Sbep[X_1^2]=\overline{\sigma}^2$ and $\cSbep[X_1^2]=\underline{\sigma}^2$. Denote $S_0^X=0$, $S_n^X=\sum_{k=1}^n X_k$, $V_0=0$, $V_n=\sum_{k=1}^n X_k^2$. And suppose that $(\widetilde{\Omega}, \widetilde{\mathscr{H}}, \widetilde{\mathbb E})$
 is a sub-linear expectation space which is rich enough such that there is a G-Brownian motion $W(t)$     with $W(1)\sim N(0,[\underline{\sigma}^2,\overline{\sigma}^2])$.  We    denote a pair of   capacities corresponding to the sub-linear expectation $\widetilde{\mathbb E}$ by  $(\widetilde{\Capc},\widetilde{\cCapc})$, and the extension of $\widetilde{\mathbb E}$ by  $\widetilde{\mathbb E}^{\ast}$.

\section{Main results}\label{sectMain}
\setcounter{equation}{0}

 We consider the convergence of the process $S_{[nt]}$. Because it is not in $C[0,1]$, it needs to be modified.   Define the $C[0,1]$-valued random variable $\widetilde{S}_n^X(\cdot)$ by setting
$$\widetilde{S}_n^X(t)=\begin{cases}
 \sum_{j=1}^k X_j,  \; \text{ if }  t=k/n \; (k=0,1,\ldots, n);\\
 \text{ extended by linear interpolation in each interval }\\
  \qquad \quad \big[[k-1]n^{-1}, kn^{-1}\big].
\end{cases}$$
Then
$ \widetilde{S}_n^X(t)=S_{[nt]}^X+(nt-[nt])X_{[nt]+1}$. Here $[nt]$ is the largest integer less than or equal to $nt$. Zhang (2015) obtained the functional central limit theorem as follows. 

 \medskip
\noindent {\bf Theorem A}  {\em Suppose $\Sbep[(X_1^2-b)^+]\to 0$ as $b\to \infty$. Then for all bounded   continuous function $\varphi:C[0,1]\to \mathbb R$,
\begin{equation}\label{eqthFCLT.1}
\Sbep\left[\varphi\left(\frac{\widetilde{S}_n^(\cdot)}{\sqrt{n}}\right)\right]\to \widetilde{\mathbb E}\left[\varphi\Big( W(\cdot) \Big) \right].
 \end{equation}
 \medskip }

Replacing the normalization factor $\sqrt{n}$ by $\sqrt{V_n}$, we obtain the self-normalized process of partial sums:
$$W_n(t)=\frac{\widetilde{S}_n^X(t)}{\sqrt{V_n}}, $$
where $\frac{0}{0}$ is defined to be $0$. 
Our  main result is   the following self-normalized  functional  central limit theorem (FCLT).
\begin{theorem}\label{thSelfFCLT} Suppose $\Sbep[(X_1^2-b)^+]\to 0$ as $b\to \infty$. Then for all bounded   continuous function $\varphi:C[0,1]\to \mathbb R$,
\begin{equation}\label{eqthSelfFCLT.1}
\Sbep^{\ast}\left[\varphi\left(W_n(\cdot)\right)\right]\to \widetilde{\mathbb E}\left[\varphi\Big(\frac{W(\cdot)}{\sqrt{\langle W  \rangle_1}}\Big) \right].
 \end{equation}
 In particular,
   for all bounded   continuous function $\varphi:\mathbb R\to \mathbb R$,
\begin{align}\label{eqthSelfFCLT.2}
\Sbep^{\ast}\left[\varphi\Big(\frac{S_n^X}{\sqrt{V_n}}\Big)\right]\to & \widetilde{\mathbb E}\left[\varphi\Big(\frac{W(1)}{\sqrt{\langle W  \rangle_1}}\Big) \right]\\
&=\sup_{\theta\in \Theta}\ep_P\left[\varphi\left(\frac{\int_0^1\theta(s) d B(s)}{\sqrt{\int_0^1 \theta^2(s) ds }}\right) \right].\nonumber
 \end{align}
\end{theorem}
\begin{remark} It is obvious that
$$\widetilde{\mathbb E}\left[\varphi\Big(\frac{W(\cdot)}{\sqrt{\langle W  \rangle_1}}\Big) \right]
\ge \ep_P \left[\varphi\Big(B(\cdot)\Big) \right]. $$
An interesting problem is how to estimate the upper bounds of   the expectations on the right hands of (\ref{eqthSelfFCLT.1}) and  (\ref{eqthSelfFCLT.2}).

Further,  $\frac{W(\cdot)}{\sqrt{\langle W\rangle_1}}\overset{d}=\frac{\overline{W}(\cdot)}{\sqrt{\langle \overline{W}\rangle_1}}$, where
$\overline{W}(t)$ is a G-Brownian motion with $\overline{W}(1)\sim N(0,[r^{-2},1])$, $r^2=\overline{\sigma}^2/\underline{\sigma}^2$.

\end{remark}

For the classical self-normalized central limit theorem, Gin\'e, G\"otze  and Mason (1997) showed that the finiteness of the second moments can be relaxed to the condition (\ref{classicalCLT4}).
Cs\"org\H o, Szyszkowicz and  Wang (2003) proved the self-normalized functional central limit theorem under (\ref{classicalCLT4}).
The next theorem gives a similar result under the sub-linear expectation and is an extension of Theorem \ref{thSelfFCLT}.
\begin{theorem} \label{thSelfFCLT2} Let $\{X_n; n\ge 1\}$ be a sequence of independent and identically distributed random variables in the sub-linear expectation space  $(\Omega, \mathscr{H}, \Sbep)$ with $\Sbep[X_1]=\cSbep[X_1]=0$. Denote
  $l(x)=\Sbep[X_1^2\wedge x^2]$. Suppose
 \begin{description}
  \item[\rm(I)]  $x^2\Capc(|X_1|\ge x)=o\big( l(x)\big)$ as $x\to \infty$;
  \item[\rm(II)] $\lim_{x\to \infty} \frac{\Sbep[X_1^2\wedge x^2]}{\cSbep[X_1^2\wedge x^2]}=r^2<\infty$;
  \item[\rm (III)] $\Sbep[(|X_1|-c)^+]\to 0$ as $c\to \infty$.
\end{description}
Then the conclusions of Theorem \ref{thSelfFCLT} remain true with $W(t)$ being a G-Brownian motion such that $W(1)\sim N(0,[r^{-2},1])$.
\end{theorem}

\begin{remark} The conditions (III) implies that $\Sbep[(|X_1|-x)^+]\le \int_x^{\infty}\Capc(|X_1|\ge y)dy$ and $=o(x^{-1}l(x))$ if the condition (I) is satisfied.
When $\Sbep$ is a continuous sub-linear expectation, then for any random variable $Y$ we have $\Sbep[|Y|]\le \int_0^{\infty}\Capc(|Y|\ge y)dy$ and so
the condition (III) can be removed.

\end{remark}

 \section{Invariance principle }\label{sectIP}
 \setcounter{equation}{0}

 To prove Theorems \ref{thSelfFCLT} and \ref{thSelfFCLT2}, we will prove a  new Donsker's invariance principle.
  Let $\{(X_i,Y_i); i\ge 1\}$ be a sequence of independent and identically distributed random vectors in the sub-linear expectation space $(\Omega, \mathscr{H}, \Sbep)$ with $\Sbep[X_1]=\Sbep[-X_1]=0$, $\Sbep[X_1^2]=\overline{\sigma}^2$, $\cSbep[X_1^2]=\underline{\sigma}^2$, $\Sbep[Y_1]=\overline{\mu}$, $\cSbep[Y_1]=\underline{\mu}$. Denote
\begin{equation}\label{eqGequation} G(p,q)=\Sbep[\frac{1}{2} q X_1^2+pY_1], \;\; p,q\in \mathbb R.
\end{equation}
Let $\xi$ be a G-normal distributed random variable, $\eta$ be  a maximal distributed  random variable such that the distribution of $(\xi,\eta)$ is   characterized by the following parabolic partial differential equation (PDE for short) defined on $[0,\infty)\times \mathbb R\times \mathbb R$:
\begin{equation}\label{eqPDE}
  \partial_t u -G\left(\partial_y u, \partial_{xx}^2 u\right) =0,
  \end{equation}
i.e., if for any bounded Lipschitz   function $\varphi(x,y):\mathbb R^2\to \mathbb R$, the function $u(x,y,t)=\widetilde{\mathbb E}\left[\varphi\left(x+\sqrt{t} \xi, y+t\eta\right)\right]$ ($x,y \in \mathbb R, t\ge 0$) is the unique viscosity solution of the  PDE (\ref{eqPDE}) with Cauchy condition $u|_{t=0}=\varphi$.

Further, let $B_t$ and $b_t$ be two random processes   such that the distribution of the process $(B_{\cdot},b_{\cdot})$ is  characterized by
  \begin{description}
  \item[\rm (i)] $B_0=0$, $b_0=0$;
    \item[\rm (ii)] for any $0\le t_1\le \ldots\le t_k\le s\le t+s$, $(B_{s+t}-B_s, b_{s+t}-b_s)$ is independent to $(B_{t_j}, b_{t_j}), j=1,\ldots,k$,   in sense that, for any $\varphi\in C_{l,Lip}(\mathbb R^{2(k+1)})$,
   \begin{align}\label{eqthFCLT.2}
    &  \widetilde{\mathbb E}\left[\varphi\left((B_{t_1}, b_{t_1}),\ldots,(B_{t_k}, b_{t_k}), (B_{s+t}-B_{s}, b_{s+t}-b_{s})\right)\right]\nonumber\\
   &\qquad  =  \widetilde{\mathbb E}\left[\psi\left((B_{t_1}, b_{t_1}),\ldots,(B_{t_k}, b_{t_k})\right)\right]
    \end{align}
    where
    $$\psi\left((x_1, y_1),\ldots,(x_k, y_k)\right)= \widetilde{\mathbb E}\left[\varphi\left((x_1, y_1),\ldots,(x_k, y_k), (B_{s+t}-B_s, b_{s+t}-b_s)\right)\right];$$
    \item[\rm (iii)]  for any $t,s>0$, $(B_{s+t}-B_s,b_{s+t}-b_s)\overset{d}\sim  (B_t,b_t)$ under $\widetilde{\mathbb E}$;
    \item[\rm (iv)]  for any $t>0$, $(B_t,b_t)\overset{d}\sim \big(\sqrt{t}B_1, tb_1\big)$ under $\widetilde{\mathbb E}$;
    \item[\rm (v) ]  the distribution of $(B_1,b_1)$ is   characterized by the PDE (\ref{eqPDE}).
 \end{description}
It is easily seen that $B_t$ is a G-Brownian motion with $B_1\sim N\big(0,[\underline{\sigma}^2,\overline{\sigma}^2]\big)$,
and $(B_t,b_t)$ is a generalized  G-Brownian motion introduced by Peng (2010a). The existence of the generalized G-Brownian motion can be found in Peng (2010a).

 \begin{theorem}\label{thFCLT}    Suppose $\Sbep[(X_1^2-b)^+]\to 0$ and $\Sbep[(|Y_1|-b)^+]\to 0$ as $b\to \infty$. Let
$$
\widetilde{\bm W}_n(t)=\left(\frac{\widetilde{S}_n^X(t)}{\sqrt{n}},  \frac{\widetilde{S}_n^Y(t)}{n}\right). $$
Then   for any  bounded continuous function $\varphi:C[0,1]\times C[0,1]\to \mathbb R$,
  \begin{equation} \label{eqCLT} \lim_{n\to \infty}\Sbep\left[\varphi\left(\widetilde{\bm W}_n(\cdot) \right)\right]= \widetilde{\mathbb E}\left[\varphi\left(B_{\cdot},b_{\cdot}\right)\right].
  \end{equation}
   Further, let   $p\ge 2$, $q\ge 1$, and assume $\Sbep[|X_1|^p]<\infty$, $\Sbep[|Y_1|^q]<\infty$. Then for all continuous function $\varphi:C[0,1]\times C[0,1]\to \mathbb R$ with $|\varphi(x,y)|\le C(1+\|x\|^p+\|y\|^q)$,
   \begin{equation} \label{eqCLT.2} \lim_{n\to \infty}\Sbep^{\ast}\left[\varphi\left(\widetilde{\bm W}_n(\cdot) \right)\right]= \widetilde{\mathbb E}\left[\varphi\left(B_{\cdot},b_{\cdot}\right)\right].
  \end{equation}
  Here $\|x\|=\sup_{0\le t\le 1}|x(t)|$ for $x\in C[0,1]$.
  \end{theorem}

 \begin{remark}\label{remark1} When $X_k$ and $Y_k$ are random vectors in $\mathbb R^d$ with $\Sbep[X_k]=\Sbep[-X_k]=0$, $\Sbep[(\|X_1\|^2-b)^+]\to 0$ and $\Sbep[(\|Y_1\|-b)^+]\to 0$
 as $b\to\infty$. Then the function $G$ in (\ref{eqGequation}) becomes
 $$ G(p,A)=\Sbep\left[\frac{1}{2}\langle AX_1,X_1\rangle+\langle p,Y_1\rangle\right],\;\; p\in \mathbb R^d, A\in\mathbb S(d), $$
 where $\mathbb S(d)$ is  the collection of all $d\times d$ symmetric matrices. The conclusion of Theorem \ref{thFCLT} remains true with
 the distribution of $(B_1,b_1)$ being characterized by the following  parabolic partial differential equation   defined on $[0,\infty)\times \mathbb R^d\times \mathbb R^d$:
$$
  \partial_t u -G\left(D_y u, D_{xx}^2 u\right) =0,\;\; u|_{t=0}=\varphi,
 $$
 where $D_y =( \partial_{y_i})_{i=1}^n$ and $D_{xx}^2=(\partial_{x_ix_j}^2)_{i,j=1}^d$.
 \end{remark}

\begin{remark} As a conclusion of Theorem \ref{thFCLT}, we have
$$ \Sbep\left[\varphi\left(\frac{S_n^X}{\sqrt{n}},\frac{S_n^Y}{n}\right)\right]\to  \widetilde{\mathbb E}\left[\varphi(B_1,b_1)\right],\;\; \varphi\in C_b(\mathbb R^2). $$
This is proved by Peng (2010a) under the conditions $\Sbep[|X_1|^{2+\delta}]<\infty$ and $\Sbep[|Y_1|^{1+\delta}]<\infty$ (c.f., Theorem 3.6 and Remark 3.8 therein).
%
\end{remark}

Before the proof, we need several lemmas. For random vectors $\bm X_n$ in $(\Omega, \mathscr{H}, \Sbep)$ and $\bm X$ in $\big(\widetilde{\Omega}, \widetilde{\mathscr{H}},\widetilde{\mathbb E}\big)$, we write $\bm X_n\overset{d}\to \bm X$ if
$$ \Sbep[\varphi(\bm X_n)]\to \widetilde{\mathbb E}[\varphi(\bm X)] $$
for any bounded continuous $\varphi$.  Write $\bm X_n \overset{\Capc}\to \bm x$  if $\Capc(\|\bm X_n-\bm x\|\ge \epsilon)\to 0$ for any $\epsilon>0$. $\{\bm X_n\}$ is called to be uniformly integrable if
$$\lim_{b\to \infty}\limsup_{n\to \infty} \Sbep[(\|\bm X_n\|-b)^+]= 0. $$
The following three lemmas are obvious.
\begin{lemma}\label{lem2} If $\bm X_n\overset{d}\to \bm X$ and $\varphi$ is a continuous function, then $\varphi(\bm X_n)\overset{d}\to \varphi(\bm X)$.
\end{lemma}

\begin{lemma}\label{lemSlutsky} {\em (Slutsky's Lemma)} Suppose $\bm X_n\overset{d}\to \bm X$, $\bm Y_n \overset{\Capc}\to \bm y$, $\eta_n\overset{\Capc}\to a$, where $a$ is a constant and $\bm y$ is a constant vector, and $\widetilde{\Capc}(\|\bm X\|>\lambda)\to 0$ as $\lambda\to \infty$. Then
$(\bm X_n, \bm Y_n, \eta_n)\overset{d}\to (\bm X,\bm y, a)$, and as a result,
$\eta_n\bm X_n+\bm Y_n\overset{d}\to a\bm X+\bm y$.
\end{lemma}
\begin{remark} Suppose $\bm X_n\overset{d}\to \bm X$. Then $\widetilde{\Capc}(\|\bm X\|>\lambda)\to 0$ as $\lambda\to \infty$ is equivalent to that $\{\bm X_n;n\ge 1\}$ is tight, i.e.,
$$ \lim_{\lambda\to \infty} \limsup_{n\to \infty} \Capc\left(\|\bm X_n\|>\lambda\right)=0, $$
because for all $\epsilon>0$, one can define a continuous function $\varphi(x)$ such that $I\{x>\lambda+\epsilon\}\le \varphi(x)\le I\{x>\lambda]$ and so
\begin{align*}
&\widetilde{\Capc}(\|\bm X\|>\lambda+\epsilon)\le    \widetilde{\mathbb E}[\varphi(\|\bm X\|)]=\lim_{n\to \infty} \Sbep[\varphi(\|\bm X_n\|)]
\le  \limsup_{n\to \infty} \Capc\left(\|\bm X_n\|> \lambda\right), \\
&\limsup_{n\to \infty} \Capc\left(\|\bm X_n\|> \lambda+\epsilon\right)\le   \lim_{n\to \infty} \Sbep[\varphi(\|\bm X_n\|)]=\widetilde{\mathbb E}[\varphi(\|\bm X\|)]
\le   \widetilde{\Capc}(\|\bm X\|> \lambda).
\end{align*}
 \end{remark}

\begin{lemma}\label{lem3} Suppose $\bm X_n\overset{d}\to \bm X$.
  \begin{description}
    \item[\rm (a)] If  $\{\bm X_n\}$ is  uniformly integrable and $\widetilde{\mathbb E}[((\|\bm X\|-b)^+]\to 0$ as $b\to \infty$,
then
\begin{equation}\label{eqlem3.1} \Sbep[\bm X_n]\to \widetilde{\mathbb E}[\bm X].
\end{equation}
    \item[\rm (b)]   If  $\sup_n\Sbep[|\bm X_n\|^q<\infty$ and  $\widetilde{\mathbb E} [|\bm X\|^q<\infty$ for some $q>1$, then (\ref{eqlem3.1}) holds.
  \end{description}
  \end{lemma}

The following lemma is proved by Zhang (2015).
\begin{lemma}\label{lem4} Suppose that $\bm X_n\overset{d}\to \bm X$, $\bm Y_n\overset{d}\to \bm Y$,
 $\bm Y_n$ is independent to $\bm X_n$ under $\Sbep$ and $\widetilde{\Capc}(\|\bm X\|>\lambda)\to 0$  and $\widetilde{\Capc}(\|\bm Y\|>\lambda)\to 0$  as $\lambda\to \infty$. Then
 $ (\bm X_n,\bm Y_n)\overset{d }\to  (\overline{\bm X},\overline{\bm Y}), $
 where $\overline{\bm X}\overset{d}=\bm X$, $\overline{\bm Y}\overset{d}=\bm Y$ and  $\overline{\bm Y}$ is independent to $\overline{\bm X}$ under $\widetilde{\mathbb E}$.
\end{lemma}

The next lemma is on the Rosenthal type inequalities  due to Zhang (2014).

\begin{lemma}\label{lemRosen}    Let $\{X_1,\ldots, X_n\}$ be a sequence  of independent random variables in $(\Omega, \mathscr{H}, \Sbep)$.
\begin{description}
  \item[\rm (a)] Suppose $p\ge 2 $. Then
\begin{align}\label{eqlemRosen1}
\Sbep\left[\max_{k\le n} \left|S_k\right|^p\right]\le  & C_p\left\{ \sum_{k=1}^n \Sbep [|X_k|^p]+\left(\sum_{k=1}^n \Sbep [|X_k|^2]\right)^{p/2} \right. \nonumber
 \\
& \qquad \left. +\left(\sum_{k=1}^n \big[\big(\cSbep [X_k]\big)^-+\big(\Sbep [X_k]\big)^+\big]\right)^{p}\right\}.
\end{align}

\item[\rm (b)] Suppose  $\Sbep[X_k]\le 0$, $k=1,\ldots, n$.    Then
\begin{equation}\label{eqlemRosen2}
\Sbep\left[\left|\max_{k\le n} (S_n-S_k)\right|^p\right]\le 2^{2-p}\sum_{k=1}^n \Sbep [|X_k|^p], \;\; \text{ for } 1\le p\le 2
\end{equation}
and
\begin{align}\label{eqlemRosen3}
\Sbep\left[\left|\max_{k\le n}(S_n- S_k)\right|^p\right]\le & C_p\left\{ \sum_{k=1}^n \Sbep [|X_k|^p]+\left(\sum_{k=1}^n \Sbep [|X_k|^2]\right)^{p/2}\right\} \nonumber\\
 \le & C_p n^{p/2-1} \sum_{k=1}^n \Sbep [|X_k|^p], \;\; \text{ for }   p\ge 2.
\end{align}
\end{description}

\end{lemma}

 \begin{lemma} \label{lemMoment}    Suppose $\Sbep[X_1]=\Sbep[-X_1]=0$ and $\Sbep[X_1^2]<\infty$.  Let $\overline{X}_{n,k}=(-\sqrt{n})\vee X_k\wedge \sqrt{n}$, $\widehat{X}_{n,k}=X_k-\overline{X}_{n,k}$, $\overline{S}_{n,k}^X=\sum_{j=1}^k \overline{X}_{n,j}$ and $\widehat{S}_{n,k}^X=\sum_{j=1}^k\widehat{X}_{n,j}$, $k=1,\ldots, n$. Then
\begin{align*}
\Sbep\left[\max_{k\le n} \left|\frac{\overline{S}_{n,k}^X}{\sqrt{n}}\right|^q\right]\le C_q, \;\; \text{ for all } q\ge 2,
\end{align*}
and
$$  \lim_{n\to \infty} \Sbep\left[\max_{k\le n} \left|\frac{\widehat{S}_{n,k}^X}{\sqrt{n}}\right|^p\right]=0 $$
whenever $\Sbep[(|X_1|^p-b)^+]\to 0$ as $b\to \infty$ if $p= 2$, and $\Sbep[|X_1|^p]<\infty$   if $p> 2$.
\end{lemma}
{\bf Proof.} Note $\Sbep[X_1]=\cSbep[X_1]=0$. So, $|\cSbep[\overline{X}_{n,1}]|=|\cSbep[X_1]-\cSbep[\overline{X}_{n,1}]|\le \Sbep|\widehat{X}_{n,1}|\le\Sbep[(|X_1|^2-n)^+]n^{-1/2}$ and
$|\Sbep[\overline{X}_{n,1}]|=|\Sbep[X_1]-\Sbep[\overline{X}_{n,1}]|\le \Sbep|\widehat{X}_{n,1}|\le \Sbep[(|X_1|^2-n)^+]n^{-1/2}$.
By Rosnethal's inequality (c.f, (\ref{eqlemRosen1})),
\begin{align*}
 & \Sbep\left[\max_{k\le n} \left|\overline{S}_{n,k}^X\right|^{q}\right]
  \le     C_p\left\{ n \Sbep [|\overline{X}_{n,1}|^q+\left(n \Sbep [|\overline{X}_{n,1}|^2]\right)^{q/2}\right.\\
    & \qquad\qquad \qquad\qquad \left.+\left(n\big[\big(\cSbep [\overline{X}_{n,1}]\big)^-+\big(\Sbep [\overline{X}_{n,1}]\big)^+\big]\right)^{q}\right\}\\
& \;\; \le  C_q\left\{ n n^{q/2-1}\Sbep [|X_1|^{2}]+n^{q/2}\left(\Sbep [X_1^2]\right)^{q/2}+\left(nn^{-1/2}\Sbep\left[(X_1^2-n)^+\right]\right)^{q}\right\} \\
& \;\; \le  C_q n^{q/2}\left\{\Sbep [|X_1|^{2}]+\left(\Sbep [X_1^2]\right)^q\right\}, \;\; \text{ for all } q\ge 2
\end{align*}
 and
\begin{align*}
  \Sbep\left[\max_{k\le n} \left|\widehat{S}_{n,k}^X\right|^{p}\right]
\le  &   C_p\left\{ n \Sbep [|\widehat{X}_{n,1}|^{p}]+\left(n \Sbep [|\widehat{X}_{n,1}|^2]\right)^{p/2}\right. \\
 & \left.\qquad +\left(n\big[\big(\cSbep [\widehat{X}_{n,1}]\big)^-+\big(\Sbep [\widehat{X}_{n,1}]\big)^+\big]\right)^{p}\right\}\\
   \le &  C_p\left\{  n\Sbep \left[\big(|X_1|^p-n^{p/2}\big)^+\right]+n^{p/2}\left(\Sbep \left[(X_1^2-n)^+\right]\right)^{p/2}\right.\\
   & \left.\qquad +n^{p/2}\left(\Sbep \left[(X_1^2-n)^+\right]\right)^{p}\right\},\; p\ge 2.
\end{align*}
The proof is completed.\hfill $\Box$

\begin{lemma} \label{lem7} (a) Suppose $p\ge 2$,  $\Sbep[X_1]=\Sbep[-X_1]=0$,  $\Sbep[(X_1^2-b)^+]\to 0$ as $b\to \infty$ and  $\Sbep[|X_1|^p]<\infty$. Then
$$ \left\{\max_{k\le n}\left|\frac{S_k^X}{\sqrt{n}}\right|^p\right\}_{n=1}^{\infty} \; \text{ is uniformly  integrable and so is tight}. $$
(b)   Suppose $p\ge 1$,   $\Sbep[(|Y_1|-b)^+]\to 0$ as $b\to \infty$, and   $\Sbep[|Y_1|^p]<\infty$. Then
$$ \left\{\max_{k\le n}\left|\frac{S_k^Y}{n}\right|^p\right\}_{n=1}^{\infty} \; \text{ is uniformly  integrable and so is tight}. $$
\end{lemma}
{\bf Proof.} (a)   follows from Lemma \ref{lemRosen}. (b) is obvious by noting
\begin{align*}
& \Sbep\left[\left(\Big(\frac{\max_{k\le n}|S_k^Y|}{n}-b\Big)^+\right)^p\right]
\le \Sbep\left[\left(\frac{\sum_{k=1}^n (|Y_k|-b)^+}{n}\right)^p\right]\\
\le& C_p \left(\frac{\sum_{k=1}^n \Sbep[(|Y_k|-b)^+]}{n}\right)^p \\
& \qquad  +C_p
\frac{\Sbep\Big[\Big|\left(\sum_{k=1}^n\{ (|Y_k|-b)^+-\Sbep[(|Y_k|-b)^+]\}\right)^+\Big|^p\Big]}{n^p}\\
 \le &  C_p\Big(\Sbep[(|Y_1|-b)^+]\Big)^p+C_p\big(n^{-p/2}+n^{1-p}\big)\Sbep\big[(|Y_k|^p-b^p)^+\big]
\end{align*}
by the Rosenthal type inequalities (\ref{eqlemRosen2}) and (\ref{eqlemRosen3}). \hfill $\Box$

\begin{lemma} \label{lem8} Suppose $\Sbep[(|Y_1|-b)^+]\to 0$ as $b\to \infty$. Then for any $\epsilon>0$,
$$\Capc\left(\frac{S_n^Y}{n}>\Sbep[Y_1]+\epsilon\right)\to 0 \text{ and } \Capc\left(\frac{S_n^Y}{n}<\cSbep[Y_1]-\epsilon\right)\to 0.  $$
\end{lemma}
{\bf Proof.} Let $Y_{k,b}=(-b)\vee Y_k\wedge b$, $S_{n,1}=\sum_{k=1}^n Y_{k,b}$ and $S_{n,2}=S_n^Y-S_{n,1}$. Note $\Sbep\big[Y_{1,b}]\to \Sbep[Y_1]$ as $b\to \infty$.
Suppose $\left|\Sbep[Y_{1,b}] - \Sbep[Y_1]\right|<\epsilon/4$. Then by Kolmogorov's  inequality (c.f. (\ref{eqlemRosen2})),
\begin{align*}
&\Capc\left(\frac{S_{n,1}}{n}>\Sbep[Y_1]+\epsilon/2\right)\le \Capc\left(\frac{S_{n,1}}{n}>\Sbep[Y_{k,b}]+\epsilon/4\right)\\
\le & \frac{16}{n^2\epsilon^2} \Sbep\left[\left(\Big(\sum_{k=1}^n \big(Y_{k,b}-\Sbep[Y_{k,b}]\big)\Big)^+\right)^2\right]
\\
\le &  \frac{32}{n^2\epsilon^2}  \sum_{k=1}^n \Sbep\left[\big(Y_{k,b}-\Sbep[Y_{k,b}]\big)^2\right]\le \frac{32(2b)^2}{n\epsilon^2}\to 0.
\end{align*}
On the other hand,
\begin{align*}
\Capc\left(\frac{S_{n,2}}{n}> \epsilon/2\right)\le & \frac{2}{n\epsilon} \sum_{k=1}^n \Sbep|Y_k-Y_{k,b}|
 \le   \frac{2}{\epsilon}\Sbep[(|Y_1|-b)^+] \to 0 \text{ as } b \to \infty.
 \end{align*}
 It follows that
 $$\Capc\left(\frac{S_n^Y}{n}>\Sbep[Y_1]+\epsilon\right)\to 0. $$
 By considering $\{-Y_k\}$ instead, we have
 $$\Capc\left(\frac{S_n^Y}{n}<\cSbep[Y_1]-\epsilon\right)=\Capc\left(\frac{-S_n^Y}{n}>\Sbep[-Y_1]+\epsilon\right)\to 0. \;\;\hfill \Box $$

 \bigskip

  {\bf Proof of Theorem \ref{thFCLT}}. We first show the tightness of $\widetilde{\bm W}_n$. It is easily seen that
  $$ w_{\delta}\left(\frac{\widetilde{S}_n^Y(\cdot)}{n}\right)
  \le 2\delta b+\frac{\sum_{k=1}^n (|Y_k|-b)^+}{n}. $$
  It follows that for any $\epsilon>0$,  if $\delta<\epsilon/(4b)$, then
  $$ \sup_n\Capc\left( w_{\delta}\left(\frac{\widetilde{S}_n^Y(\cdot)}{n}\right)\ge \epsilon\right)
  \le \sup_n \Capc\left(\sum_{k=1}^n (|Y_k|-b)^+\ge n\frac{\epsilon}{2}\right)\le \frac{2}{\epsilon}\Sbep\left[(|Y_1|-b)^+\right]. $$
Letting $\delta\to 0$ and then $b\to \infty$ yields
$$\sup_n\Capc\left( w_{\delta}\left(\frac{\widetilde{S}_n^Y(\cdot)}{n}\right)\ge \epsilon\right)\to 0 \text{ as } \delta \to 0.$$
For any $\eta>0$, we choose $\delta_k\downarrow 0$ such that, if
 $$A_k=\left\{x: \omega_{\delta_k}(x)<\frac{1}{k}\right\}, $$
then $\sup_n\Capc\left(  \widetilde{S}_n^Y(\cdot)/n \in A_k^c\right)\le \eta/2^{k+1}$.
Let $A=\{x:|x(0)|\le a\}$, $K_2=A\bigcap_{k=1}^{\infty}A_k$. Then by the Arzel\'a-Ascoli theorem, $K_2\subset C_b(C[0,1])$  is compact.  It is obvious that
$\{ \widetilde{S}_n^Y(\cdot)/n\not\in A\}=\emptyset$ since $ \widetilde{S}_n^Y(0)/n=0$. Next, we show that
$$\Capc( \widetilde{S}_n^Y(\cdot)/n\in K^c)\le \Capc( \widetilde{S}_n^Y(\cdot)/n\in A^c)+\sum_{k=1}^{\infty}\Capc( \widetilde{S}_n^Y(\cdot)/n\in A_k^c). $$
 Note that when $\delta<1/(2n)$,
$$ \omega_{\delta}( \widetilde{S}_n^Y(\cdot)/n)\le 2n |t-s|\max_{i\le n}  |Y_i|/n \le 2   \delta  \max_{i\le n} |Y_i|. $$
Choose a $k_0$ such that $\delta_k<1/(2Mk)$ for $k\ge k_0$. Then on the event $ E=\{ \max_{i\le n} |Y_i|\le M\}$, $\{ \widetilde{S}_n^Y(\cdot)/n\in A_k^c\}=\emptyset$ for $k\ge k_0$. So,
by the (finite) sub-additivity of $\Capc$,
\begin{align*}
&\Capc\big(E \bigcap \{ \widetilde{S}_n^Y(\cdot)/n\in K^c\}\big)\\
\le & \Capc\big(E \bigcap\{ \widetilde{S}_n^Y(\cdot)/n \in  A^c\}\big)+\sum_{k=1}^{k_0}\Capc\big(E\bigcap \{ \widetilde{S}_n^Y(\cdot)/n \in A_k^c\}\big) \\
\le & \Capc( \widetilde{S}_n^Y(\cdot)/n \in A^c)+\sum_{k=1}^{\infty}\Capc( \widetilde{S}_n^Y(\cdot)/n \in A_k^c).
\end{align*}
On the other hand,
$$\Capc(E^c )\le \frac{\Sbep[\max_{i\le n} |Y_i|]}{M}\le \frac{n\Sbep[|Y_1|]}{M}. $$
It follows that
\begin{align*}
\Capc\big(  \widetilde{S}_n^Y(\cdot)/n\in K_2^c \big)
\le  \Capc( \widetilde{S}_n^Y(\cdot)/n \in A^c)+\sum_{k=1}^{\infty}\Capc( \widetilde{S}_n^Y(\cdot)/n \in A_k^c)+ \frac{n \Sbep[|Y_1|]}{M}.
\end{align*}
Letting $M\to \infty$ yields
\begin{align*}  \Capc\big( \widetilde{S}_n^Y(\cdot)/n\in K_2^c \big)\le &
   \Capc( \widetilde{S}_n^Y(\cdot)/n \in A^c)+\sum_{k=1}^{\infty}\Capc( \widetilde{S}_n^Y(\cdot)/n \in A_k^c)\\
< &  0+\sum_{k=1}^{\infty} \frac{\eta}{2^{k+1}}<\frac{\eta}{2}.
\end{align*}
We conclude that for any $\eta>0$, there exists a  compact $K_2\subset C_b(C[0,1])$ such that
 \begin{equation}\label{eqProofTight.1}
  \sup_n \Sbep^{\ast}\left[I\left\{\frac{\widetilde{S}_n^Y(\cdot)}{n}\not\in K_2\right\}\right]=\sup_n \Capc\left \{\frac{\widetilde{S}_n^Y(\cdot)}{n}\not\in K_2\right\}<\eta/2.
  \end{equation}
 Next, we show that for any $\eta>0$, there exists a compact $K_1\subset C_b(C[0,1])$ such that
 \begin{equation}\label{eqProofTight.2} \sup_n \Sbep^{\ast}\left[I\left\{\frac{\widetilde{S}_n^X(\cdot)}{\sqrt{n}}\not\in K_1\right\}\right]=\sup_n \Capc\left \{\frac{\widetilde{S}_n^X(\cdot)}{\sqrt{n}}\not\in K_1\right\}<\eta/2.
 \end{equation}
  Similar to (\ref{eqProofTight.1}), it is sufficient to show that
\begin{equation}\label{eqProofTight.3} \sup_n\Capc\left( w_{\delta}\left(\frac{\widetilde{S}_n^X(\cdot)}{\sqrt{n}}\right)\ge \epsilon\right)\to 0 \text{ as } \delta \to 0.
\end{equation}
   With the same argument of Billingsley (1968, Pages 56-59, c.f., (8.12)), for large $n$,
   \begin{align*}
   &\Capc\left( w_{\delta}\left(\frac{\widetilde{S}_n^X(\cdot)}{\sqrt{n}}\right)\ge 3\epsilon\right)
     \le   \frac{2}{\delta} \Capc\left(\max_{i\le [n\delta]} \frac{|S_i^X|}{\sqrt{[n\delta]}}\ge  \epsilon  \frac{\sqrt{n}}{\sqrt{[n\delta]} }  \right) \\
     \le &\frac{2}{\delta} \Capc\left(\max_{i\le [n\delta]} \frac{|S_i^X|}{\sqrt{[n\delta]}}\ge \frac{\epsilon }{\sqrt{2 \delta }}  \right) \le  \frac{4}{\epsilon^2}\Sbep\left[\left(\max_{i\le [n\delta]} \Big|\frac{S_i^X}{\sqrt{[n\delta]}}\Big|^2-\frac{\epsilon^2 }{ 2 \delta  }\right)^+\right].
  \end{align*}
  It follows that
  $$ \lim_{\delta\to 0} \limsup_{n\to \infty} \Capc\left( w_{\delta}\left(\frac{\widetilde{S}_n^X(\cdot)}{\sqrt{n}}\right)\ge 3\epsilon\right)=0$$
  by Lemma \ref{lem7} (a) where $p=2$.
  On the other hand, for fixed $n$,  if $\delta<1/(2n)$ then
$$ \omega_{\delta}( \widetilde{S}_n^Y(\cdot)/n)\le 2n |t-s|\max_{i\le n}  |X_i|/\sqrt{n} \le 2   \delta \sqrt{n} \max_{i\le n} |X_i|. $$
It follows that
$$ \lim_{\delta\to 0}   \Capc\left( w_{\delta}\left(\frac{\widetilde{S}_n^X(\cdot)}{\sqrt{n}}\right)\ge  \epsilon\right)=0  $$
for each $n$. It follows that (\ref{eqProofTight.3}) holds.

 Now, by combing (\ref{eqProofTight.1}) and (\ref{eqProofTight.2}) we obtain the tightness of $\widetilde{\bm W}_n$ as follows.
 \begin{equation} \label{eqTightness}\sup_n \Sbep^{\ast}\left[I\left\{\widetilde{\bm W}_n(\cdot)\not\in K_1\times K_2\right\}\right]<\eta.
 \end{equation}
 Define $\Sbep_n$ by
 $$ \Sbep_n[\varphi]=\Sbep\left[\varphi\big(\widetilde{W}_n(\cdot)\big)\right],\;\; \varphi\in C_b\big(C[0,1]\times C[0,1]\big). $$
 Then the sequence of sub-linear expectations $\{\Sbep_n\}_{n=1}^{\infty}$ is tight by (\ref{eqTightness}). By Theorem 9 of Peng (2010b),  $\{\Sbep_n\}_{n=1}^{\infty}$ is weakly compact, namely, for each subsequence $\{\Sbep_{n_k}\}_{k=1}^{\infty}$, $n_k\to \infty$, there
exists a further subsequence  $\{\Sbep_{m_j}\}_{j=1}^{\infty} \subset  \{\Sbep_{n_k}\}_{k=1}^{\infty}$, $m_j\to \infty$,  such that, for each $\varphi\in C_b\big(C[0,1]\times C[0,1]\big)$,
 $\{\Sbep_{m_j}[\varphi]\}$ is
a Cauchy sequence. Define ${\mathbb F}[\cdot]$ by
$$ {\mathbb F}[\varphi]=\lim_{j\to \infty}\Sbep_{m_j}[\varphi], \; \varphi\in C_b\big(C[0,1]\times C[0,1]\big). $$
Let $\overline{\Omega}=C[0,1]\times C[0,1]$, and $(\xi_t,\eta_t)$ be the the canonical process $\xi_t(\omega) = \omega_t^{(1)}$, $\eta_t(\omega)=\omega_t^{(2)}$  $(\omega=(\omega^{(1)},\omega^{(2)})\in \widetilde{\Omega})$. Then
\begin{equation}\label{eqSubsequence}\Sbep\left[\varphi\big(\widetilde{W}_{m_j}(\cdot)\big)\right]\to {\mathbb F}\left[\varphi(\xi_{\cdot},\eta_{\cdot})\right],\;\;\varphi\in C_b\big(C[0,1]\times C[0,1]\big).
\end{equation}
The topological completion of $C_b(\overline{\Omega})$  under the Banach norm ${\mathbb F}[\|\cdot\|]$ is denoted by $L_{\mathbb F} (\overline{\Omega})$. ${\mathbb F}[\cdot]$ can be extended uniquely to
a sublinear expectation on  $L_{\mathbb F} (\overline{\Omega})$.

Next, it is sufficient to show that $(\xi_t,\eta_t)$ defined on the sub-linear space $(\overline{\Omega}, L_{\mathbb F} (\overline{\Omega}), {\mathbb F})$ satisfies (i)-(v) and so $(\xi_{\cdot},\eta_{\cdot})\overset{d}=(B_{\cdot},b_{\cdot})$, which means that the limit distribution of any subsequence of $\widetilde{\bm W}_n(\cdot)$ is uniquely  determined.

(i) is obvious.

Let $0\le t_1\le \ldots\le t_k\le s\le t+s$. By (\ref{eqSubsequence}), for any bounded continuous function $\varphi:\mathbb R^{2(k+1)}\to \mathbb R$ we have
\begin{align*}
& \Sbep\left[\varphi\big(\widetilde{W}_{m_j}(t_1),\ldots,  \widetilde{W}_{m_j}(t_k), \widetilde{W}_{m_j}(s+t)-\widetilde{W}_{m_j}(s)\big)\right]  \\
\to &{\mathbb F} \left[\varphi\big( (\xi_{t_1},\eta_{t_1}), \ldots, (\xi_{t_k},\eta_{t_k}),(\xi_{s+t}-\xi_s,\eta_{s+t}-\eta_s)\big)\right].
 \end{align*}
 Note
\begin{align*} & \sup_{0\le t\le 1}\frac{|\widetilde{S}_n^X(t)-S_{[nt]}^X|}{\sqrt{n}}\le \frac{\max_{k\le n}|X_k|}{\sqrt{n}}\overset{\Capc}\to 0,\\
 &\sup_{0\le t\le 1}\frac{|\widetilde{S}_n^Y(t)-S_{[nt]}^Y|}{n}\le \frac{\max_{k\le n}|Y_k|}{n}\overset{\Capc}\to 0.
\end{align*}
It follows that by Lemmas \ref{lemSlutsky} and \ref{lem7},
\begin{align}\label{eqfinite}
& \Sbep\left[\varphi\left(\Big(\frac{S_{[m_jt_1]}^X}{\sqrt{m_j}},\frac{S_{[m_jt_1]}^Y}{m_j}\Big),\ldots, \Big(\frac{S_{[m_jt_k]}^X}{\sqrt{m_j}},\frac{S_{[m_jt_k]}^Y}{m_j}\Big), \Big(\frac{S_{[m_j(s+t)]}^X-S_{[m_js]}^X}{\sqrt{m_j}},\frac{S_{[m_j(s+t)]}^Y-S_{[m_js]}^Y}{m_j}\Big)\right)\right]\nonumber \\
& \qquad \to  {\mathbb F} \left[\varphi\big( (\xi_{t_1},\eta_{t_1}), \ldots, (\xi_{t_k},\eta_{t_k}),(\xi_{s+t}-\xi_s,\eta_{s+t}-\eta_s)\big)\right].
 \end{align}
 In particular,
\begin{align*}  \left(\frac{S_{[m_j(s+t)]-[m_js]}^X}{\sqrt{m_j}},\frac{S_{[m_j(s+t)]-[m_js]}^Y}{m_j}\right) \overset{d}=& \left(\frac{S_{[m_j(s+t)]}^X-S_{[m_js]}^X}{\sqrt{m_j}}, \frac{S_{[m_j(s+t)]}^Y-S_{[m_js]}^Y}{m_j}\right)\\
&\overset{d}\to \big(\xi_{s+t}-\xi_s, \eta_{s+t}-\eta_s\big).
 \end{align*}
 It follows that
\begin{equation}\label{eqCLTproof} \left(\frac{S_{[m_jt]}^X}{\sqrt{m_j}},\frac{S_{[m_jt]}^Y}{m_j}\right) \overset{d}\to \big(\xi_{s+t}-\xi_s, \eta_{s+t}-\eta_s\big).
\end{equation}
 Hence,
\begin{equation}\label{eqEqualDistribution} {\mathbb F}\left[\phi(\xi_{s+t}-\xi_s,\eta_{s+t}-\eta_s)\right]={\mathbb F}\left[\phi(\xi_t,\eta_t)\right]\;\; \text{ for all } \phi\in C_b(\mathbb R^2).
\end{equation}

 Next, we show that
\begin{equation}\label{eqmoment} {\mathbb F}[|\xi_{s+t}-\xi_s|^p]\le C_p t^{p/2}\; \text{ and }\; {\mathbb F}[|\eta_{s+t}-\eta_s|^p]\le C_p t^p,\;\;\text{ for all } p\ge 2 \text{ and } t,s \ge 0.
\end{equation}
  By Lemma \ref{lem8},
   \begin{equation}\label{eqboundeta}
   \widetilde{\cCapc}\big(t\underline{\mu}-\epsilon\le \eta_{s+t}-\eta_s\le t\overline{\mu}+\epsilon\big)=1\;\; \text{ for all } \; \epsilon>0.
   \end{equation}
    It follows that
$$ {\mathbb F}[|\eta_{s+t}-\eta_s|^p]\le t^p\big|\Sbep[|Y_1|]\big|^p.
$$

  On the other hand, let $\overline{S}_{n,k}^X$ and $\widehat{S}_{n,k}^X$ be defined as in Lemma \ref{lemMoment}. Then $S_k^X=\overline{S}_{n,k}^X+
  \widehat{S}_{n,k}^X$. By (\ref{eqCLTproof}) and Lemmas \ref{lemMoment} and \ref{lemSlutsky},
$$  \frac{\overline{S}_{[m_jt],[m_jt]}^X}{\sqrt{m_j}}  \overset{d}\to  \xi_{s+t}-\xi_s\; \text{ and }\;
\Sbep\left[\left|\frac{\overline{S}_{[m_jt],[m_jt]}^X}{\sqrt{m_j}}\right|^{p}\right]\le C_p t^{p/2},\; p\ge 2. $$
It follows that
$$  {\mathbb F}\left[| \xi_{s+t}-\xi_s|^p\wedge b\right]=\lim_{n\to\infty}\Sbep\left[\left|\frac{\overline{S}_{[m_jt],[m_jt]}^X}{\sqrt{m_j}}\right|^{p}\wedge b\right]\le C_p t^{p/2}, \;\text{   for any } b>0. $$
Hence,
$$ {\mathbb F}\left[| \xi_{s+t}-\xi_s|^p \right]=\lim_{b\to \infty}  {\mathbb F}\left[| \xi_{s+t}-\xi_s|^p\wedge b\right]\le C_p t^{p/2} $$ by the completeness of $(\overline{\Omega}, L_{\mathbb F} (\overline{\Omega}), {\mathbb F})$. (\ref{eqmoment}) is proved.

  Now,  note that $(X_i,Y_i), i=1,2,\ldots$, are independent and identically distributed. By (\ref{eqfinite}) and Lemma \ref{lem4},
 it is easily seen that $(\xi_{\cdot}, \eta_{\cdot})$ satisfies (\ref{eqthFCLT.2})  for   $\varphi\in C_b(\mathbb R^{2(k+1)})$.
 Note that, by (\ref{eqmoment}), the random variables concerned in (\ref{eqthFCLT.2}) and (\ref{eqEqualDistribution})  have finite moments of each order. The function space $C_b(\mathbb R^{2(k+1)})$ and $C_b(\mathbb R^2)$ can be extended to $C_{l,Lip}(\mathbb R^{2(k+1)})$ and
 $C_{l,Lip}(\mathbb R^2)$, respectively, by  elemental arguments. So, (ii) and (iii) is proved.

 For (iv) and (v), we let $\varphi:\mathbb R^2\to \mathbb R$ be a bounded Lipschitz function and consider
 $$ u(x,y,t)={\mathbb F}\left[\varphi(x+\xi_t,y+\eta_t)\right]. $$
 It is sufficient to show that $u$ is a viscosity solution of  the PDE (\ref{eqPDE}). In fact, due to the uniqueness of the viscosity solution, we will have
 $${\mathbb F}\left[\varphi(x+\xi_t,y+\eta_t)\right]=\widetilde{\mathbb E}\left[\varphi(x+\sqrt{t} \xi,y+t\eta)\right], \;\; \varphi\in C_{b,Lip}(\mathbb R^2). $$
 Taking $x=0$ and  $y=0$ yields (iv) and (v).

 To verify   the PDE (\ref{eqPDE}), firstly it is easily seen that
 $$ \Sbep\left[\frac{q}{2}   \Big(\frac{S_{[nt]}^X}{\sqrt{n}}\Big)^2+p \frac{S_{[nt]}^Y}{n}\right]=\frac{[nt]}{n}\Sbep\left[\frac{q}{2}   \Big(\frac{S_{[nt]}^X}{\sqrt{[nt]}}\Big)^2+p \frac{S_{[nt]}^Y}{[nt]}\right]=\frac{[nt]}{n}G(p,q).$$
 Note that $\left\{\frac{q}{2}   \Big(\frac{S_{[nt]}^X}{\sqrt{n}}\Big)^2+p \frac{S_{[nt]}^Y}{n}\right\}$ is uniformly integrable by Lemma \ref{lem7}. By Lemma \ref{lem3} we conclude that
 $${\mathbb F}\left[\frac{q}{2}\xi_t^2+p\eta_t\right]=\lim_{m_j\to \infty}\Sbep\left[\frac{q}{2}   \Big(\frac{S_{[m_jt]}^X}{\sqrt{m_j}}\Big)^2+p \frac{S_{[m_jt]}^Y}{m_j}\right]=t G(p,q). $$
 Also, it is easy to verify that
 $ |u(x,y,t)-u(\overline{x},\overline{y},t)|\le C ( |x-\overline{x}|+|y-\overline{y}|) $,
  $ |u(x,y,t)-u(x,y,s)|\le C\sqrt{|t-s|}$ by Lipschitz continuity of $\varphi$, and
\begin{align*}
  u(x,y,t)=&{\mathbb F}\left[\varphi(x+\xi_s+\xi_t-\xi_s,y+\eta_s+\eta_t-\eta_s)\right]\\
  =& {\mathbb F}\left[ {\mathbb F}
  \left[\varphi(x+\overline{x}+\xi_t-\xi_s,y+\overline{y}+\eta_t-\eta_s)\right]\big|_{(\overline{x},\overline{y})=(\xi_s,\eta_s)}\right]\\
  = &
  {\mathbb F}\left[u(x+\xi_s,y+\eta_s, t-s)\right],\; 0\le s\le t.
  \end{align*}
  Let $\psi(\cdot,\cdot,\cdot)\in C_b^{3,3,2}(\mathbb R,\mathbb R,[0,1])$  be a smooth function with $\psi\ge u$ and $\psi(x,y,t)=u(x,y,t)$. Then
  \begin{align*}
  0=& {\mathbb F}\left[u(x+\xi_s,y+\eta_s, t-s)-u(x,y,t)\right]\le  {\mathbb F}\left[\psi(x+\xi_s,y+\eta_s, t-s)-\psi(x,y,t)\right]\\
  = &  {\mathbb F}\left[\partial_x\psi(x,y,t)\xi_s+\frac{1}{2} \partial_{xx}^2\psi(x,y,t)\xi_s^2+\partial_y\psi(x,y,t)\eta_s-\partial_t \psi(x,y,t) s+I_s\right]\\
  \le &    {\mathbb F}\left[\partial_x\psi(x,y,t)\xi_s+\frac{1}{2} \partial_{xx}^2\psi(x,y,t)\xi_s^2+\partial_y\psi(x,y,t)\eta_s-\partial_t \psi(x,y,t) s\right]+ {\mathbb F}[|I_s|]\\
  =&  {\mathbb F}\left[\frac{1}{2} \partial_{xx}^2\psi(x,y,t)\xi_s^2+\partial_y\psi(x,y,t)\eta_s\right]-\partial_t \psi(x,y,t) s+ {\mathbb F}[|I_s|]\\
  =& sG(\partial_y\psi(x,y,t),\partial_{xx}^2(x,y,t)) -s\partial_t \psi(x,y,t)+ {\mathbb F}[|I_s|],
  \end{align*}
  where
  $$|I_s|\le C\big( |\xi_s|^3+|\eta_s|^2+s^2\big). $$
  By (\ref{eqmoment}), we have
  $  {\mathbb F}[|I_s|]\le C\big( s^{3/2}+s^2+s^2\big)=o(s). $
  It follows that $[\partial_t \psi- G(\partial_y\psi,\partial_{xx}^2)](x,y,t)\le 0$.  Thus $u$ is a viscosity subsolution of (\ref{eqPDE}). Similarly we can prove that $u$ is a viscosity supersolution of (\ref{eqPDE}). Hence (\ref{eqCLT}) is proved.

As for (\ref{eqCLT.2}), let $\varphi:C[0,1]\times C[0,1]\to \mathbb R$ be a continuous function with $|\varphi(x,y)|\le C_0(1+\|x\|^p+\|y\|^q)$. For $\lambda>4 C_0$, let $\varphi_{\lambda}(x,y)=(-\lambda)\vee (\varphi(x,y)\wedge \lambda)\in C_b(C[0,1])$. It is easily seen that $\varphi(x,y)=\varphi_{\lambda}(x,y)$ if $|\varphi(x,y)|\le \lambda$.   If $|\varphi(x,y)|>\lambda$, then
\begin{align*}
   |\varphi(x,y)- & \varphi_{\lambda}(x,y)|=|\varphi(x,y)|-\lambda\le  C_0(1+\|x\|^p+\|y\|^q)-\lambda\\
 \le & C_0\Big\{\Big(\|x\|^p-\lambda/(4C_0)\Big)^+ +\Big(\|y\|^q-\lambda/(4C_0)\Big)^+\Big\}.
 \end{align*}
Hence
$$|\varphi(x,y)-\varphi_{\lambda}(x,y)|
 \le   C_0\Big\{\Big(\|x\|^p-\lambda/(4C_0)\Big)^+ +\Big(\|y\|^q-\lambda/(4C_0)\Big)^+\Big\}.$$
 It follows that
\begin{align*}
&\lim_{\lambda\to \infty}  \limsup_{n\to \infty}\left|\Sbep^{\ast}\left[\varphi\left(\widetilde{\bm W}_n(\cdot) \right)\right]-
 \Sbep\left[\varphi_{\lambda}\left(\widetilde{\bm W}_n(\cdot) \right)\right]\right| \\
\le &  \lim_{\lambda\to \infty}  \limsup_{n\to \infty} C_0\left\{ \Sbep\left[\left(\max_{k\le n}\left|\frac{S_k^X}{\sqrt{n}}\right|^p-\frac{\lambda}{4C_0}\right)^+\right]+\Sbep\left[\left(\max_{k\le n}\left|\frac{S_k^Y}{n}\right|^q-\frac{\lambda}{4C_0}\right)^+\right]\right\}\\
=&0,
\end{align*}
by Lemma \ref{lem7}. On the other hand, by (\ref{eqCLT}),
$$\lim_{n\to \infty}\Sbep\left[\varphi_{\lambda} \left(\widetilde{\bm W}_n(\cdot) \right)\right]= \widetilde{\mathbb E}\left[\varphi_{\lambda} \left(B_{\cdot},b_{\cdot}\right)\right]\to \widetilde{\mathbb E}\left[\varphi \left(B_{\cdot},b_{\cdot}\right)\right]\;\; \text{ as } \lambda\to \infty. $$
  (\ref{eqCLT.2}) is proved, and the proof Theorem \ref{thFCLT} is now completed. \hfill $\Box$

 \bigskip
   {\bf Proof of Remark \ref{remark1}.} When $X_k$ and $Y_k$ are $d$-dimensional random vectors, the tightness (\ref{eqTightness}) of $\widetilde{\bm W_n}(\cdot)$ also follows because each sequence of the components  of the vector $\widetilde{\bm W_n}(\cdot)$ is tight. Also, (\ref{eqmoment}) remains true because each component has this property. On the other hand, it follows that
   \begin{align*}
{\mathbb F}\left[\frac{1}{2}\left\langle A\xi_t,\xi_t\right\rangle+\left\langle p,\eta_t\right\rangle\right]
   = & \lim_{m_j\to \infty} \Sbep\left[\frac{1}{2}\left\langle A\frac{S_{[m_jt]}^X}{\sqrt{m_j}},\frac{S_{[m_jt]}^X}{\sqrt{m_j}}\right\rangle+\left\langle p,\frac{S_{[m_jt]}^Y}{m_j}\right\rangle\right]\\
    = & \lim_{m_j\to \infty} \frac{[m_jt]}{m_j} G(p,A)=t  G(p,A).
   \end{align*}
   The remainder proof is the same as that of Theorem \ref{thFCLT}. \hfill $\Box$

 \section{Proof of the self-normalized FCLTs}\label{sectProof}
 \setcounter{equation}{0}
  Let $Y_k=X_k^2$. The function $G(p,q)$ in (\ref{eqGequation}) becomes
 $$  G(p,q)=\Sbep\left[\Big(\frac{q}{2} +p\Big) X_1^2\right]=\Big(\frac{q}{2} +p\Big)^+\overline{\sigma}^2-\Big(\frac{q}{2} +p\Big)^-\underline{\sigma}^2, \;\; p,q\in \mathbb R.
$$
 Then the process $(B_t, b_t)$ in (\ref{eqCLT}) and the process $(W(t), \langle W\rangle_t)$  are identically distributed.

  In fact, note
 $$  \langle W\rangle_{t+s}-\langle W\rangle_t=(W(t+s)-W(t))^2-2\int_0^s (W(t+x)-W(t))d(W(t+x)-W(t)). $$
 It is easy to verify that  $(W(t), \langle W\rangle_t)$ satisfies (i)-(iv) for $(B_{\cdot}, b_{\cdot})$. It remains to show that
 $(B_1, b_1)\overset{d}= (W(1), \langle W\rangle_1)$.
 Let $\{X_n;n\ge 1\}$ be a sequence of independent and identically distributed random variables with $X_1\overset{d}= W(1)$. Then by Theorem \ref{thFCLT},
 $$ \left(\frac{\sum_{k=1}^nX_k}{\sqrt{n}},\frac{\sum_{k=1}^n X_k^2}{n}\right)\overset{d}\to (B_1, b_1). $$
 On the other hand, let $t_k=\frac{k}{n}$. Then
$$ \left(\frac{\sum_{k=1}^nX_k}{\sqrt{n}},\frac{\sum_{k=1}^n X_k^2}{n}\right)\overset{d}=
  \left(W(1), \sum_{k=1}^n (W(t_k)-W(t_{k-1}))^2\right)\overset{L_2}\to (W(1), \langle W\rangle_1). $$
  Hence $(B_{\cdot},b_{\cdot})\overset{d}=(W(\cdot), \langle W\rangle_{\cdot})$. We conclude the following proposition from Theorem \ref{thFCLT}.
  \begin{proposition} \label{prop1} Suppose   $\Sbep[(X_1^2-b)^+]\to 0$ as $b\to \infty$. Then  for all bounded   continuous function
  $\psi:C[0,1]\times C[0,1]\to \mathbb R$,
 $$
\Sbep\left[\psi\left(\frac{\widetilde{S}_n^X(\cdot)}{\sqrt{n}},\frac{\widetilde{V}_n (\cdot)}{n}\right)\right]\to \widetilde{\mathbb E}\left[\psi\Big( W(\cdot), \langle W  \rangle_{\cdot} \Big) \right]
$$
 where $\widetilde{V}_n(t)=V_{[nt]}+(nt-[nt])X^2_{[nt]+1}$,
 and in particular, for all bounded   continuous function
  $\psi:C[0,1]\times \mathbb R\to \mathbb R$,
\begin{equation}\label{eqprop1.1}
\Sbep\left[\psi\left(\frac{\widetilde{S}_n^X(\cdot)}{\sqrt{n}},\frac{V_n}{n}\right)\right]\to \widetilde{\mathbb E}\left[\psi\Big( W(\cdot), \langle W  \rangle_1 \Big) \right].
 \end{equation}
   \end{proposition}
  \bigskip

  Now, we begin the proof of Theorem \ref{thSelfFCLT}. Let $a=\underline{\sigma}^2/2$ and $b=2\overline{\sigma}^2$.  According to  (\ref{eqboundeta}), we have $\widetilde{\cCapc}\big( \underline{\sigma}^2-\epsilon <
  \langle W  \rangle_1<\overline{\sigma}^2+\epsilon\big)=1$ for all $\epsilon>0$. Let $\varphi:C[0,1]\to \mathbb R$ be a bounded continuous function. Define
  $$ \psi\big(x(\cdot),y\big)=\varphi\left(\frac{x(\cdot)}{\sqrt{ a\vee y\wedge b}}\right), \;\; x(\cdot)\in C[0,1],\;y\in \mathbb R. $$
  Then $\psi:C[0,1]\times \mathbb R\to \mathbb R$ is a bounded continuous function. Hence by Proposition \ref{prop1},
  $$ \Sbep\left[\varphi\left(\frac{\widetilde{S}_n^X(\cdot)/ \sqrt{n}}{\sqrt{a\vee(V_n/n)\wedge b}}\right)\right]\to \widetilde{\mathbb E}\left[\varphi\left(\frac{W(\cdot)}{\sqrt{a\vee (\langle W  \rangle_1)\wedge b}} \right) \right]=\widetilde{\mathbb E}\left[\varphi\left(\frac{W(\cdot)}{\sqrt{   \langle W  \rangle_1)} } \right) \right]. $$
  On the other hand,
\begin{align*}
  \limsup_{n\to \infty}   & \left|\Sbep^{\ast}\left[\varphi\left(\frac{\widetilde{S}_n^X(\cdot)/ \sqrt{n}}{\sqrt{  V_n/n }}\right)\right]-\Sbep\left[\varphi\left(\frac{\widetilde{S}_n^X(\cdot)/ \sqrt{n}}{\sqrt{a\vee(V_n/n)\wedge b}}\right)\right]\right|\\
  \le & C \limsup_{n\to \infty} \Capc\left(V_n/n\not\in (a,b)\right)\\ \le & C\widetilde{\Capc}\left(
  \langle W  \rangle_1\ge 3\overline{\sigma}^2/2\right) + C\widetilde{\Capc}\left(
  \langle W  \rangle_1\le 2\underline{\sigma}^2/3\right)=0.
 \end{align*}
 It follows that
   $$ \Sbep^{\ast}\left[\varphi\left(\frac{\widetilde{S}_n^X(\cdot)}{\sqrt{ V_n }}\right)\right]\to  \widetilde{\mathbb E}\left[\varphi\left(\frac{W(\cdot)}{\sqrt{   \langle W  \rangle_1)} } \right) \right]. $$
   The proof is now completed.\hfill $\Box$

\bigskip
{\bf Proof of Theorem \ref{thSelfFCLT2}}.  Firstly,  note that
\begin{align*} \Sbep[X_1^2\wedge x^2]\le & \Sbep[X_1^2\wedge(kx)^2]\le \Sbep[X_1^2\wedge x^2]+k^2x^2\Capc(|X_1|>x),\;\; k\ge 1,  \\
 \Sbep[|X_1|^r\wedge x^r]\le &\Sbep[|X_1|^r\wedge (\delta x)^r]+\Sbep[(\delta x)^r\vee |X_1|^r\wedge x^r]\\
\le &  \delta^{r-2} x^{r-2}l(\delta x)+x^r \Capc(|X_1|\ge \delta x), \;\;0<\delta<1,\; r>2.
\end{align*}
The condition (I) implies that $l(x)$ is slowly varying as $x\to \infty$ and
$$ \Sbep[|X_1|^r\wedge x^r]=o(x^{r-2}l(x)), \; r>2. $$
Further
$$ \frac{\Sbep^{\ast}[X_1^2I\{|X_1|\le x\}]}{l(x)}\to 1, $$
$$ C_{\Capc}\big(|X_1|^rI\{|X_1|\ge x\}\big)=\int_{x^r}^{\infty} \Capc(|X_1|^r\ge y)dy =o(x^{2-r} l(x)),\;\; 0<r<2.$$
If the conditions (I) and (III) are satisfied, then
$$ \Sbep[(|X_1|-x)^+]\le \Sbep^{\ast}[|X_1|I\{|X|\ge x\}] \le C_{\Capc}\big(|X_1|I\{|X_1|\ge x\}\big)=o(x^{-1} l(x)). $$

Now, let $d_t=\inf\{x: x^{-2}l(x)=t^{-1}\}$. Then $nl(d_n)=d_n^2$. Similar to Theorem \ref{thSelfFCLT}, it is sufficient to show that for all bounded   continuous function
  $\psi:C[0,1]\times C[0,1]\to \mathbb R$,
$$ \Sbep\left[\psi\left(\frac{\widetilde{S}_n^X(\cdot)}{d_n},\frac{\widetilde{V}_n(\cdot)}{d_n^2}\right)\right]\to \widetilde{\mathbb E}\left[\psi(W(\cdot), \langle W\rangle_{\cdot})\right]\;\; \text{ with }  W(1)\sim N(0,[r^{-2}, 1]). $$
Let $\overline{X}_k=\overline{X}_{k,n}=(-d_n)\vee X_k\wedge d_n$, $\overline{S}_k =\sum_{i=1}^k \overline{X}_i$, $\overline{V}_k=\sum_{i=1}^k  \overline{X}_i^2$. Denote $ \overline{S}_n(t)=\overline{S}_{[nt]}+(nt-[nt])\overline{X}_{[nt]+1}$ and $\overline{V}_n(t)=\overline{V}_{[nt]}+(nt-[nt])\overline{X}^2_{[nt]+1}$. Note
$$\Capc\left(X_k\ne \overline{X}_k \text{ for some } k\le n\right)\le n \Capc\left(|X_1|\ge d_n\right)=n\cdot o\Big(\frac{l(d_n)}{d_n^2}\Big)=o(1). $$
It is sufficient to show that for all bounded   continuous function
  $\psi:C[0,1]\times C[0,1]\to \mathbb R$,
$$ \Sbep\left[\psi\left(\frac{\overline{S}_n(\cdot)}{d_n},\frac{\overline{V}_n(\cdot)}{d_n^2}\right)\right]\to \widetilde{\mathbb E}\left[\psi(W(\cdot), \langle W\rangle_{\cdot})\right]. $$
Following the line of the proof of Theorem \ref{thFCLT}, we need only to show that
\begin{description}
  \item[\rm (a)] For any $0<t\le 1$,
  $$\limsup_{n\to \infty} \Sbep\left[\max_{k\le [nt]}\left|\frac{\overline{S}_k}{d_n}\right|^p\right]\le C_p t^{p/2},\;\;   \limsup_{n\to \infty}\Sbep\left[\max_{k\le [nt]}\left|\frac{\overline{V}_k}{d_n^2}\right|^p\right]\le C_p t^p,\;\; \forall p\ge 2;$$
  \item[\rm (b)]  For any $0<t\le 1$,
 $$ \lim_{n\to \infty} \Sbep\left[ \frac{q}{2}   \Big(\frac{\overline{S}_{[nt]}}{d_n}\Big)^2+p \frac{\overline{V}_{[nt]}}{d_n^2}\right]=tG(p,q), $$
 where
 $$G(p,q)=\Big(\frac{q}{2}+p\Big)^+ -  r^{-2}\Big(\frac{q}{2}+p\Big)^-; $$
  \item[\rm (c)]
  $$ \max_{k\le n} \frac{|X_k|}{d_n}\overset{\Capc}\to 0. $$
\end{description}
In fact, (a)   implies the tightness of $\left(\frac{\widetilde{S}_n^X(\cdot)}{d_n},\frac{\widetilde{V}_n(\cdot)}{d_n^2}\right)$ and (\ref{eqmoment}), and (b) implies the distribution of the limit process is uniquely determined.

Firstly, (c) is obvious since
$$ \Capc\Big(\max_{k\le n}|X_k|\ge \epsilon d_n\Big)\le n \Capc\Big( |X_1|\ge \epsilon d_n\Big)=o(1) n \frac{l(\epsilon d_n)}{\epsilon^2 d_n^2}
=o(1) n \frac{l( d_n)}{d_n^2}=o(1). $$

As for (a), by the Rosenthal type inequality (\ref{eqlemRosen1}),
\begin{align}\label{eqproofthad.5}
&\Sbep  \left[\max_{k\le [nt]}\left|\frac{\overline{S}_k}{d_n}\right|^p\right]
\le  C_pd_n^{-p}\left\{ [nt] \Sbep\left[|X_1|^p\wedge d_n^p\right]+\Big( [nt] \Sbep\left[|X_1|^2\wedge d_n^2\right]\Big)^{p/2}\right.\nonumber\\
 & \qquad \qquad \qquad + \left. \Big([nt] (\cSbep[(-d_n)\vee X_1\wedge d_n])^++[nt] (\Sbep[(-d_n)\vee X_1\wedge d_n])^+\Big)^p\right\}\nonumber\\
 \le& C_pd_n^{-p}\left\{ [nt] \Sbep\left[|X_1|^p\wedge d_n^p\right]+\Big( [nt] \Sbep\left[|X_1|^2\wedge d_n^2\right]\Big)^{p/2} + \Big([nt] \Sbep[(|X_1|-d_n)^+] \Big)^p\right\}\nonumber\\
  \le& C_pd_n^{-p}\left\{ [nt] o\Big(d_n^{p-2}l(d_n)\Big) +\Big( [nt] l(d_n)\Big)^{p/2} + \bigg([nt] o\Big(\frac{l(d_n)}{d_n}\Big) \bigg)^p\right\}\nonumber\\
 =& o(1)[nt]\frac{l(d_n)}{d_n^2}+\Big(\frac{[nt]}{n}
 \Big)^{p/2} \Big(\frac{nl(d_n)}{d_n^2}\Big)^{p/2}+o(1)\Big( [nt] \frac{l(d_n)}{d_n^2}\Big)^p\le C_pt^{p/2}+o(1),
\end{align}
and similarly,
\begin{align*}
 \Sbep  \left[\max_{k\le [nt]}\left|\frac{\overline{V}_k}{d_n^2}\right|^p\right]
\le & C_pd_n^{-2p}\left\{ [nt] \Sbep\left[|X_1|^{2p}\wedge d_n^{2p}\right]+\Big( [nt] \Sbep\left[|X_1|^4\wedge d_n^4\right]\Big)^{p/2}\right.\\
 & + \left. \Big([nt]  \cSbep[  X_1^2\wedge d_n^2] +[nt] (\Sbep[  X_1^2\wedge d_n^2] \Big)^p\right\}\\
  =& o(1)+ C_p \Big( [nt] \frac{l(d_n)}{d_n^2}\Big)^p\le C_pt^{p}+o(1).
\end{align*}
Thus (a) follows.

As for (b), note
\begin{align*}
 \frac{q}{2}   \Big(\frac{\overline{S}_{[nt]}}{d_n}\Big)^2+p \frac{\overline{V}_{[nt]}}{d_n^2}
 =\Big(\frac{q}{2}+p\Big)\frac{\overline{V}_{[nt]}}{d_n^2}+q\frac{\sum_{k=1}^{[nt]-1}\overline{S}_{k-1}\overline{X}_k}{d_n^2}.
\end{align*}
By (\ref{eqproofthad.5}),
\begin{align*}
&\Sbep\left[ \sum_{k=1}^{[nt]-1}\overline{S}_{k-1}\overline{X}_k \right]\le  \sum_{k=1}^{[nt]-1}\Sbep\left[\overline{S}_{k-1}\overline{X}_k \right]\\
\le &  \sum_{k=1}^{[nt]-1}\left\{\Sbep\left[( \overline{S}_{k-1})^+\right]\Sbep[\overline{X}_k]-\Sbep\left[( \overline{S}_{k-1})^-\right]\cSbep[\overline{X}_k]\right\}\\
\le & \sum_{k=1}^{[nt]-1} \left(\Sbep\left[| \overline{S}_{k-1}|^2\right]\right)^{1/2}\Sbep[(|X_1|-d_n)^+]\\
=& O\Big((nl(d_n)^{1/2}\Big)\cdot n\Sbep[(|X_1|-d_n)^+]\\
=& O(d_n)\cdot n\cdot o\Big(\frac{l(d_n)}{d_n}\Big)=o(d_n^2),
\end{align*}
and similarly,
$$ \Sbep\left[- \sum_{k=1}^{[nt]-1}\overline{S}_{k-1}\overline{X}_k \right]=o(d_n^2). $$
On the other hand,
$$ \frac{\Sbep[V_{[nt]}]}{d_n^2}=\frac{[nt]\Sbep[X_1^2\wedge d_n^2]}{d_n^2}=\frac{[nt]}{n}\frac{nl(d_n)}{d_n^2}=\frac{[nt]}{n}\to t $$
and
$$ \frac{\cSbep[V_{[nt]}]}{d_n^2}= \frac{[nt]\cSbep[X_1^2\wedge d_n^2]}{d_n^2}=\frac{[nt]}{n}\frac{\cSbep[X_1^2\wedge d_n^2] }{\Sbep[X_1^2\wedge d_n^2]}\to t r^{-2}. $$
Hence we conclude that
\begin{align*}
  \Sbep& \left[\frac{q}{2}   \Big(\frac{\overline{S}_{[nt]}}{d_n}\Big)^2+p \frac{\overline{V}_{[nt]}}{d_n^2}\right]
 =  \Sbep\left[\Big(\frac{q}{2}+p\Big)\frac{\overline{V}_{[nt]}}{d_n^2}\right]+o(1)\\
&\quad =  t  \left[ \Big(\frac{q}{2}+p\Big)^+ -  r^{-2}\Big(\frac{q}{2}+p\Big)^-\right] +o(1).
\end{align*}
Thus (b) is verified, and the proof is completed. \hfill $\Box$

\end{document}